# The Mathematisation of the World

## Uncovering the Socio-Economic Tensions for Ethics in Mathematics Education

Dennis Müller[1]

**Abstract:**

The mathematisation of the socio-economic sphere, where mathematics actively constructs social reality, presents a challenge for studies on ethics in mathematics and its education. While existing scholarship on ethics in mathematics offers insights, it often remains philosophically driven and disconnected from other relevant disciplines. This paper addresses this gap by asking how debates on ethics in mathematics and its education can be connected with economic sociology, and what socio-economic tensions become visible through this connection. Drawing from concepts such as imagined futures, varieties of capitalism, and variegated capitalism, we synthesise a new perspective. This analysis reveals six interconnected tensions: a socio-economic valuation gap regarding ethics education; the multifaceted implementation of mathematics across different capitalist systems; its material opaqueness; a growing gap between economic power and social unaccountability; the enclosure of imagination limiting sustainable futures; and the erosion of multilateralism, which challenges critical pedagogy. The paper's contribution is a first step towards a structural socio-economic framework that links the limited literature on ethics in mathematics with these broader sociological perspectives.

**Keywords:** Ethics in Mathematics, Sociology of Mathematics, Critical Mathematics Education, Sustainable Mathematics Education, Economic Sociology

---

[1] Institute of Mathematics Education, University of Cologne, Germany. dennis.mueller@uni-koeln.de

# Introduction

> *"The undoubted benefits and intrinsic virtues of mathematics should not blind us to potential collateral damage as the huge juggernaut of mathematics rolls over education and society towards a reshaped future. Because of its great power and influence across the whole of society we need to undertake an ethical audit."*
>
> *Paul Ernest (2019a, p. 1).*

The "collateral damage" of mathematics that Ernest (2019a) describes is not a future problem; it is a present reality, embedded in global socio-economic structures. But why does it so often persist, go unnoticed, or get dismissed? What are the deeper socio-economic structures enabling this damage to happen? We argue that to grasp the deeper economic roots behind many of today's pressing concerns of ethics in mathematics, we must first grasp its sociology. By grounding the critique of this paper in economic sociology, we wish to describe how mathematics does not merely describe our socio-economic sphere of life but actively constructs it. We aim to present a deep entanglement of abstract mathematics with socio-economic concerns that quickly creates winners and losers.

Several tensions related to the implementation of sustainability into the mathematical curriculum can be identified (Meyer, forthcoming[2]): different modes of argumentation between mathematics and sustainability, differences between clarity and ambiguity, tensions of focus and curriculum planning, challenges in the framing of contexts, and aspects of learner polarisation. We wish to complement this practitioner-focused analysis by highlighting tensions that arise for mathematics education as a result of today's mathematisation of the socio-economic sphere. We see a two-fold need for this: first, the existing landscape of studies regarding ethical and sustainable mathematics is becoming increasingly philosophically complex, so that a new programmatic argument can provide a navigational aid; something to agree with and critique, and thus to locate one's own position in relation to (see also Müller et al., 2025). Second, despite its increasing complexity, the existing scholarship on ethics in mathematics is very recent (Ernest, 2024), and what exists is often very philosophically driven (Müller, 2025b). Therefore, it could benefit from stronger integration with insights from other subjects.

---

[2] These tensions are from a forthcoming paper that is discussed in Müller et al. (2025).

The literature of this "ethical turn"[3] appears dominated by three primary currents: philosophical inquiry, critical pedagogy (including non-Western concerns), and practical concerns. None of these, however, constitutes an application of economic sociology (cf. Ernest, 2024; Müller, 2025b), which seeks to analytically explain the institutional and social structures of economic life. This paper aims to overcome this gap and complement the more philosophical outlooks on ethics in mathematics and its education by incorporating explicit positions from economic sociology: most notably concepts from varieties of capitalism, variegated capitalism, imagined futures, and institutional isomorphisms.

Indeed, typical (philosophically and historically inspired) stories on ethics (and social justice) in mathematics tell a narrative along the following lines: absolutist philosophical positions hold that mathematical knowledge is always certain, universal, and objective. This stance unites diverse schools of thought, including Platonism, Logicism, and Formalism, around the fundamental belief that mathematical truths represent discoveries rather than human creations, and such a stance would need to be scrutinised by more socio-critical scholarship (Ernest, 1991a, 1991b). Through European colonialism, the philosophical idea of mathematical absolutism, and universalism in particular, was transformed into a concrete historical force. The colonial period witnessed the global spread of a particular mathematical approach, rooted in ancient Greek rationality, formal arithmetic, and logical proof. Mathematics became an essential instrument of imperial control, enabling colonial administration, facilitating resource exploitation, and restructuring societies along new economic lines (Bishop, 1990). This expansion created a worldwide intellectual hierarchy that systematically marginalised other forms of mathematical practice (cf. Hodgkin, 2005; Joseph, 2011; Bishop, 1990).

For much of the 20th century, most of the history of mathematics was written through a Eurocentric perspective, maintaining stories of "linear" growth, neglecting vast mathematical achievements from other cultures, and occasionally also ignoring that even Western mathematics had its periods of doubt and fear, for example, as evidenced in the 19th-century foundational crisis (Gray, 2004). In particular, over the last decades, cultural histories of mathematics have emerged, which put the spotlight on achievements of non-Western cultures and the variety of practices found in different mathematical communities today and in the past (Rowe & Dauben, 2024). The study of ethnomathematics (D'Ambrosio, 2016), the shifting focus of the philosophy of mathematics from ontological questions to mathematical practice (Ernest, 1991b; Wagner, 2017), and the interest in a socially-just sociology of mathematics (Barany & Kremakova, 2020) are part of this change in how mathematics is increasingly understood as a form of human action.

---

[3] This "ethical turn" is to be differentiated from previous scholarship on social justice in mathematics education. This social justice-centred literature has, for instance, made much use of Marxist or Foucauldian lenses. The literature on the ethical turn, in contrast, is more concerned with the ethical responsibilities of individual mathematicians, students, and educators (e.g., Ernest, 2024; Müller, 2024), a conversation which this paper seeks to ground in a structural, socio-economic framework.

Inspired by such narratives, the paper asks the following research question: How can we connect existing debates on ethics in mathematics and its education with economic sociology, and what socio-economic tensions become visible when we do so?

By investigating how economic phenomena are socially constructed and institutionally embedded, economic sociology challenges an abstract view of the economy that presents it as a separate sphere, fully governed by rationality. Instead, it seeks to examine the complex interplay of markets, institutions, culture, and power structures. Perhaps most important for this study, the field also scrutinises how tools like mathematisation are not neutral descriptors but active forces in organising economic life and legitimising specific futures. This perspective allows us to move beyond philosophical abstraction and to trace some of the socio-economic roots of ethical challenges in mathematics and its education.

As economic sociology is itself a broad and fragmented field, and not a single doctrine (cf. Maurer, 2021), we will have to engage with several of its different strands to at least attempt to be somewhat comprehensive. But even our connections mostly trace back to a few selected traditions within the field. We fully acknowledge that each of the connections identified warrants its own, deeper investigation in future work. This paper is explicitly positioned as an "exploratory" study, and a continuation of the work by Müller et al. (2025) on the triangle of ethical concerns related to mathematics. Its primary aim is to test the value of forging a connection between ethics in mathematics education and the field of economic sociology, seeking to determine what can be gained from this synthesis.

More abstractly, this paper attempts an unconventional argumentative path. As illustrated in Figure 1, we seek to connect mathematical concerns (the nature of the subject) and community concerns (pedagogy, practice) by first traversing large-scale socio-planetary concerns (the economic structures of capitalism). This approach deliberately inverts the trajectory usually taken by scholars in mathematics education, who often begin with the immediate, local concerns of the community (such as the classroom or practitioners) and then expand outwards to broader societal issues, or with a combination of both concerns (cf. Müller et al., 2025). By starting instead with the macro-level socio-economic structures related to mathematisation, this paper can be understood as a proof of concept for the analytical power of arguing in this alternative direction.

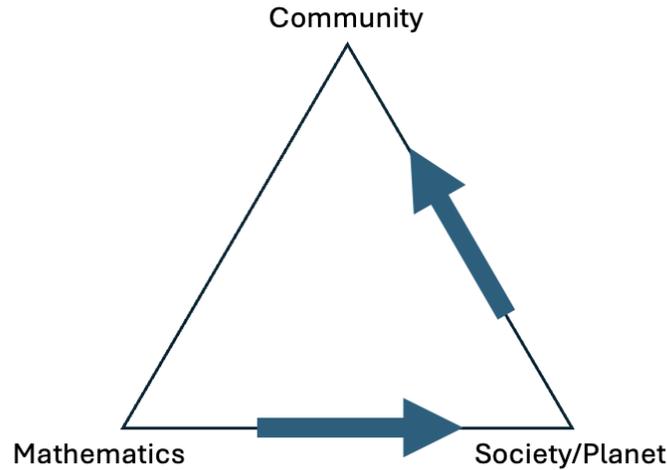

Figure 1: The trajectory of this paper[4]

To answer our research question, we proceed in three parts; the latter two of which are split into smaller steps. First, we explain our usage of the term "mathematisation" (part I), which we distinguish from quantification — the production and usage of numbers/data. Then we make the connection between economic sociology and mathematising the world (part II). Here, we specify how processes of mathematisation are related to different varieties of capitalism and variegated forms of capitalism, putting particular emphasis on the Global North/South divide. Having established this foundation, we will outline some consequences of it by looking at aspects of power, human identity, and the recent decline in multilateralism (part III). In particular, this perspective will ultimately lead us to consider how absolutist and socio-critical stances have different implied models of governance. We conclude by arguing that it may be beneficial to take on a critical pragmatic position to deal with the ethical and sustainability challenges arising from mathematics.

For each section, we consider its impact on mathematics education and the role of mathematical practitioners, outlining a tension between a social and economic focus. This will lead us to describe six socio-economic tensions. These tensions manifest in different forms — some as market dynamics, others as structural gaps, or geopolitical pressures — but they all converge to create fundamental, unresolved conflicts for the pedagogy and practice of mathematics:

1. The socio-economic valuation gap
2. The multifaceted implementation of mathematics
3. The material opaqueness of mathematics
4. The economic power versus social unaccountability gap
5. The enclosure of imagination

---

[4] This triangle is inspired by Müller et al. (2025).

6. The erosion of multilateralism undermines the critical challenger

To visualise these tensions, we will use the "Ethical and Sustainable Concerns Triangle" (ESCT) developed by Müller et al. (2025). This triangle depicts the three prototypical ethical concerns of mathematicians and educators about mathematics, the community, and the socio-planetary through a triangle to map discursive patterns, compare educational positions, and to explain different tensions in the field. Mathematical concerns are understood to include concerns about the integrity and continuity of mathematical knowledge, community concerns are, among others, social and ethical questions within and about the community of mathematicians, educators, and students, and socio-planetary concerns are about the impact of and role that mathematics plays in the socio-planetary sphere at large.

# Part I: Mathematisation versus quantification

> *"Most more advanced mathematical concepts, such as complex numbers, algebras, linear operators, Borel sets [...] were so devised that they are apt subjects on which the mathematician can demonstrate his ingenuity and sense of formal beauty. In fact, the definition of these concepts, with a realization that interesting and ingenious considerations could be applied to them, is the first demonstration of the ingeniousness of the mathematician who defines them."*
>
> *Eugene P. Wigner (1960, p. 3)*

To build our argument, we must first draw a clear distinction between quantification and the broader concept we term "mathematisation." The sociology of quantification presents a rich diversity of perspectives that predominantly centre around four questions (Berman & Hirschman, 2018, p. 258)

1. "What shapes the production of numbers?"
2. "When and how do numbers matter?"
3. "How do we govern quantification? How *should* we govern quantification?"
4. "How should scholars *study* quantification?"

A foundational work in the sociology of quantification by Espeland and Stevens (2008, p. 401) defines quantification as "the production and communication of numbers". Their framing of quantification is inspired by speech acts, building on Austin's (1975) philosophy of speech — most notably, his philosophy of locutionary, illocutionary, and perlocutionary acts. Numbers can then be locutionary acts ("saying something") as statements of fact, as descriptions, and as assertions. Here, mathematics and its various subfields (including statistics) present the "grammar" and "vocabulary" for using numbers. Illocutionary acts focus on "doing something". The illocutionary acts of numbers involve attention, persuasion, and categorisation. Perlocutionary acts then describe the aftermath of, or the effects of, using numbers: What does using numbers do to the world?

At the heart of quantification are acts of commensuration, i.e., the transformation of diverse qualities into shared metrics. Through this, otherwise incomparable aspects can be captured. Espeland and Stevens (2008) further identify five sociological dimensions through which quantification exerts its transformative influence: the labour it demands, its reactivity (the way measurement alters what is measured), its disciplinary function, its foundation for authority, and its aesthetic qualities. They conclude that quantification itself reconfigures the reality it wants to describe. Ernest (2019a, p. 11) expands on this relationship by noting that "[m]athematics is not just the measure of pre-existing reality, it is constitutive in the making of a new and enacted reality. Mathematics and measurement are entangled in and complicit in the remaking of reality we experience and inhabit."

We follow Ernest's observations and argue that, for our analysis, quantification is vital, but to properly capture the deep-seated influence of mathematical thought on the socio-economic sphere and to then draw lessons about the ethics of mathematics and its education, we think it is too narrow. Mennicken and Espeland (2019) note that the dominant intellectual traditions of quantification studies are focused on statistics (e.g., the French school), focused on the consequences of statistics and their means of production (e.g., history and philosophy of science), and on overcoming neutrality myths (e.g., critical accounting studies). Despite the field's recent expansion to include algorithms, big data, artificial intelligence, and datafication (Berman & Hirschman, 2018; Camargo & Daniel, 2021; Chun & Son, 2025; Diaz & Didier, 2016; Mennicken & Espeland, 2019; Saltelli & Di Fiore, 2020), the account of its foundational traditions confirms a strong conceptual focus on numbers rather than on mathematics.

However, because our primary audience consists of mathematicians and mathematics educators, we feel a definition focused on "the production and communication of numbers" (Espeland & Stevens, 2008, p. 401) or which has to ask, "Are algorithms quantification?" (Berman & Hirschman, 2018, p. 258), is insufficient. To have a meaningful ethical discussion with the community of mathematicians and educators, we must use a more encompassing concept that captures the impact of the entire discipline — including its underlying logics, pure and abstract concepts, and professional values and its boundaries — not just its most visible numerical outputs.[5] We therefore distinguish between quantification and the broader logic of mathematisation.

We see quantification as the most visible and widespread form of mathematisation. It is the process of turning the world into numbers. However, mathematisation, as used here, denotes the comprehensive embrace of mathematical practices and logics to structure reality. This perspective allows us to analyse not only the numbers themselves (quantification) but also the underlying rationales that are imported into the socio-economic sphere. In this framework, mathematisation is not a neutral tool; it is the ideological procedure that reframes complex social phenomena into something that can be quantified, optimised, predicted, and, thus, ultimately controlled. While quantification is number-centric, our concept of mathematisation is logic-centric. We are interested in how the adoption of mathematical rationales — such as abstraction, formal proof, and optimisation — shapes the socio-economic world, even when a specific quantification is not the end result. Algorithmic governance, for instance, is a form of mathematisation not just because it uses numbers, but because it embeds a formal, logical, and often optimising process at the heart of decision-making (cf. Issar & Aneesh, 2022).

To translate this into educational terms: a focus on quantification might critique how standardised testing scores affect school funding, whereas a focus on mathematisation critiques the underlying logic that complex human learning can be reduced to mathematical terms in the first place. This critical usage of mathematisation must be distinguished from its specific pedagogical and curricular meaning. In education, mathematisation typically describes the student's process of translating a real-world problem into mathematical terms (horizontal mathematisation) or operating within the mathematical domain to solve it (vertical mathematisation) (Jablonka & Gellert, 2007; building on Treffers, 1987).

This process, however, represents only a one-way translation from the real world to the world of mathematics.[6] The concept we are concerned with is broader, encompassing a crucial second aspect: how that very model is embedded with politics, imbues institutions with its instrumental logic, and ultimately affects the very world it sought to model. This reflexive, world-altering dimension of mathematisation is generally absent from the curriculum. For instance, the review by Makramalla et al. (2025, p. 540) builds on a systematic review by Vásquez et al. (2023) to note that there are "very few initiatives" related to the critical integration of sustainability. Textbook problems are already fixed when the students meet them; they do not change simply because a student solves them. The students, however, may of course be affected by the act of solving, as Ernest (2019a) argues when he speaks about the collateral damage of learning mathematics. This circular perspective on mathematisation is also noticed by Davis and Hersh (1986). They see the process and aftermath of the abstraction and formalisation of the real world, and wonder if it constitutes a "social tyranny".

---

[5] Similarly, Saltelli and Di Fiore (2020) note that, as the sociology of quantification has become more diverse, the need to move to an ethics of quantification requires a more encompassing conception of quantification.

[6] Jablonka & Gellert (2007, p. 6), similarly argue that, "[m]athematics is a means for the generation of new realities not only by providing descriptions of 'real world situations', but also by colonising, permeating, and transforming reality. Models become the reality, which they set out to model. Consequently, any discussion of mathematisation has to take into account the social process by which mathematical models are developed, implemented, accepted, and obscured."

Furthermore, a narrow focus on quantification, which lends itself to applied mathematics, risks missing how the values derived from pure mathematics also shape the socio-economic sphere. The esteem given to pure abstraction, the belief in a single optimal solution, or the authority granted to formal methods are all mathematical values. These values and logics are embedded in decision-making and become part of the mathematisation of the world — sometimes, entirely independent of a connection to a specific number or processes of quantification. This transformation of mathematical thinking in abstractions into "realised abstractions" that become part of our reality creates a form of "implicit mathematics" (Keitel et al., 1993; Jablonka & Gellert, 2007). In this state, the underlying mathematical abstractions are fixed in technologies and socio-economic structures, often becoming invisible and no longer reflected upon (Chevallard, 2007).

Within the literature on ethics in mathematics, there are strong debates about whether mathematics is merely a neutral tool that is applied, and thus if all ethical analysis should be restricted to its user, or whether mathematics, as a social and cultural phenomenon, and Mathematics, as an abstract Platonist realm of reality, both have embedded ethics and politics within them. A mere restriction to quantification would not capture many of these aspects, such as how mathematicians can create or discover (depending on one's ontological perspective on this matter) new mathematics that appears completely irrelevant to physical reality. We, however, think that it is important to be as encompassing as possible, for two reasons: 1) (applied) mathematics is moving extremely quickly these days, and we do not know what happens and becomes relevant in the future, and 2) mathematicians have a tendency to erect (artificial) boundaries between themselves and those areas of mathematics that come under ethical scrutiny.[7] Speaking of ethics in mathematics education, rather than ethics in applied mathematics education, requires us to be open to these perspectives, too.

---

[7] For example, see the discussion on whether mathematics counts as a profession in Müller et al. (2022), as well as general discussions on the boundaries of science in Gieryn (1999).

Doing so will also help to illuminate that the power of mathematisation does not come from any individual innovation; rather, the revolutionary impact of mathematics is driven by deeper logics and rationales. Nonetheless, the sociology of quantification will provide useful micro-level insights to the macro-level perspectives opened up by economic sociology. It is not an argument against the sociology of quantification; rather that for our purposes we need an extension of it.

The values inherent in mathematical practice itself become organising principles for society. For instance, the high value placed on axiomatic rigour and formal proof — hallmarks of pure mathematics — can translate into rigid bureaucratic structures. When organisations prioritise formalised procedures over contextual understanding, they adopt a mathematical logic that privileges abstraction over social complexity. This reliance on idealised models, even when they poorly fit the situation, can contribute to the "tyranny of metrics" (Muller, 2018) and enable the "accountability sinks" (D. Davis, 2024) we will discuss later.

Consider Ord's (2025) work on "evaluating the infinite" as an illustration of the difference. Ord employs hyperreal numbers to assign values to infinite utility streams in economics and ethics, aiming to resolve paradoxes that emerge when standard mathematics treats different infinite summations as equivalent. His approach extends beyond simple quantification as it represents a sophisticated form of mathematisation. By incorporating the mathematical ontology of hyperreals, Ord establishes new foundations for decision theory and ethics. These foundations do not just establish new forms of mathematical logic, but also of thought and rationality. This shift cannot be fully captured by classical sociology of quantification, which primarily examines the social functions of finite numbers. However, Ord's work still reinforces Espeland and Stevens's (2008) call that a sociology of quantification requires an accompanying ethics of quantification. This need becomes even more pressing with such forms of mathematisations: if one decides to adopt Ord's hyperreal framework, it would embed significant ethical and political commitments about the future into decision-making processes while asserting that infinite future utilities can be rationally addressed without discounting the future (a contestable claim). While numbers are at the centre of Ord's logic, his shift to a different number system presents a substantial shift in mathematisation that goes beyond what is ordinarily discussed when one speaks of quantification.

Having established our core analytical lens, we can proceed to apply it. The mathematisation of the socio-economic sphere does not happen randomly, nor does it happen in a vacuum. Therefore, we must ask: what are the large-scale socio-economic structures that deploy this logic, and for what purpose do they do it? To answer this question, we will now connect it to some foundational perspectives from economic sociology on (economic) decision making and the structure of capitalism(s) itself.

# Part II: The socio-economic nature of mathematising today's world

> *"Physics, because of its astonishing success at predicting the future behavior of material objects from their present state, has inspired most financial modeling. Physicists study the world by repeating the same experiments over and over again [...] It's a different story with finance and economics, which are concerned with the mental world of monetary value."*
>
> <div align="right">*Emanuel Derman & Paul Wilmott (2009)*</div>

The 21st century is witnessing a transformation in how societies are organising themselves: higher mathematics is increasingly shaping the public sphere, mediating civic life, and governing individual opportunities (O'Neil, 2016). People increasingly inhabit what Airoldi (2022) calls a "machine habitus". This mathematical governance has benefited from specific cultural shifts in modern societies, and can build on a long tradition of what Desrosières (1998) calls the "politics of large numbers". Historically, religious narratives provided the framework for both public order and private morality. However, religion has become an increasingly private matter and modern societies separate rationality and spirituality/religiosity (Taylor, 2018; Williams, 2015). Coinciding with this shift, and without a unifying framework able to capture both private and public spheres, societies increasingly turn to algorithms and other forms of mathematics to manage complexity (Marichal, 2025). The development of statistics, in particular, is closely linked to the development of public administration and state power (Desrosières, 1998).

The increase in governance by mathematics can further be traced to specific institutional pressures. DiMaggio and Powell (1983) introduced the concept of institutional isomorphisms to analyse these dynamics. Although their theory was originally developed to explain homogenisation within specific organisational fields (e.g., why universities in a field resemble one another), its logic may be plausibly extended to the macro-societal embrace of mathematisation itself. The three pressures they identify are all visible in the push for mathematisation. Coercive pressures (forces compelling conformity) exist through mandates like standardised international testing (e.g., PISA), which compel nations to align their curricula (Ernest, 2019a), and through the integration of national economies with global financial markets (Pilbeam, 2023). Mimetic pressures (copying perceived success) occur when institutions adopt strategies like data-driven management practices often disseminated by consulting firms (cf. Bogdanich & Forsythe, 2022; Marktanner, 2023; McKenna, 2010). Finally, normative pressures (shared professional values) arise from academic discourse and an overarching neoliberal ideology that equates mathematisation with modernisation, efficiency, and rationality (Amadae, 2016; Ernest, 2019a; Marktanner, 2023). "[E]conomization, calculation, measurement and valuation" (Davies, 2017, p. 22 in Ernest, 2019a) are core components of neoliberalism, driving mathematics into a direction that distances it from human experience (Ernest, 2019a).

The means of mathematisation, including models, forecasts, statistics, and algorithms, are not neutral instruments for revealing an objective reality; rather, Beckert (2016) argues they function as powerful "cognitive devices" and forms of storytelling that make particular "fictional expectations" about the future appear plausible, credible, and legitimate, even when the mathematics cannot fully capture the underlying complex reality. In Foucauldian (2020) terms, mathematisation thus operates as "regimes of truth", where power is exercised not through force, but by defining the mechanisms and discourses, including their logic, that separate true and credible statements from false ones. In essence, these mathematical tools are not just attempting to predict the future; they are creating the narratives that motivate present action under conditions of uncertainty. For example, Marktanner (2023, p. 15) demonstrates how German public authorities strategically deployed management consultants to "depoliticise" conflicts, using their reports (which were grounded in the logics of efficiency and optimisation) as "argumentative ammunition" to frame politically motivated reforms as objective necessities. Another example is given by the staggering valuations of AI startups. These are rarely based on current profits, but on an imagined future of their potential for market dominance, and, in particular, on mathematisations of their existing successes, such as scaling laws which suggest that more data and bigger models lead to better results (cf. Floridi, 2024).

As such, forms of mathematisation build the basis of how actors make decisions in 21st-century capitalist societies (Beckert, 2016; Beckert & Bronk, 2019). Crucially, Beckert distinguishes this concept — where fictions enable action despite uncalculable uncertainty — from performativity. Many modes of mathematisation possess performative qualities beyond shaping decisions. For example, when a critical mass of actors adopts a model, their collective actions can make the world behave as the model predicts, leading to systemic feedback loops (Chiodo & Müller, 2025). When mathematisation is deployed as a cognitive device, it inherently privileges that which is measurable and marginalises qualitative knowledge and tacit understanding. Muller (2018) thus speaks of a tyranny of metrics to describe how limiting these narratives can become, as the resulting "metric fixation" means that people and institutions begin to focus on the metric itself, rather than on the underlying goal and why it should be pursued. In particular, framing complex challenges to be solved through mathematics risks neglecting valuable other perspectives, including a long list of technomoral values necessary to navigate civic life (cf. Vallor, 2016).

These intense, socio-economically motivated processes of mathematisation act as irritations on the educational and mathematical (research) socio-communicative systems (Müller et al., 2025). As the demand for mathematics education that serves utilitarian goals increases, it privileges those involved in training students for competitive markets (e.g., as is evidenced by the hiring of mathematically-trained experts in machine learning). Additionally, the marginalisation of non-quantifiable objectives strengthens the normative need for more progressive and public educators, focused on teaching aspects necessary for civic life (e.g., statistical literacy). However, as these may be seen as a challenge to the existing socio-economic status quo, they may not experience the same economic appreciation.

Müller et al. (2025) build on Ernest (1991a) to map different educator archetypes onto the ESCT. As economic forces exert a strong utilitarian pull, they favour the roles of industrial trainers and technological pragmatists who service market needs and drive quantifiable efficiencies, as well as those who maintain and provide high-level mathematical research and teaching, such as the old humanist mathematicians and academic trainers (highlighted in Figure 2).

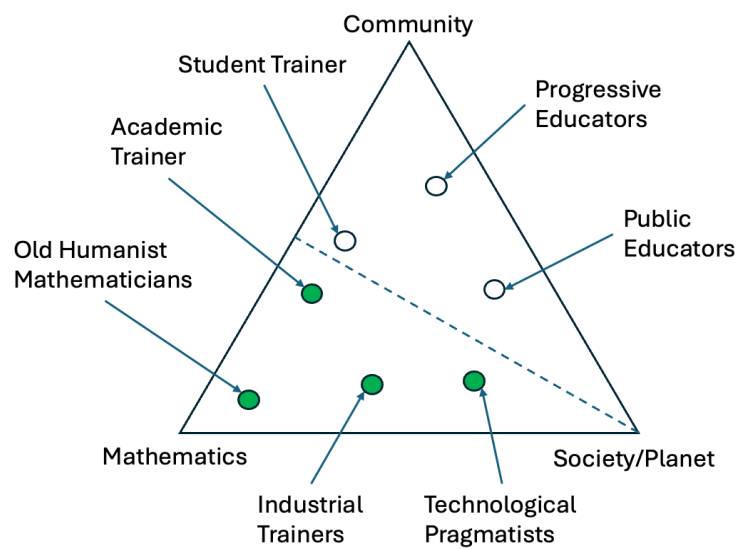

Figure 2: Educator archetypes benefitting from economically-driven mathematisation

Simultaneously, this shift creates a complex dynamic for roles oriented towards social justice that actively critique mathematics, such as progressive and public educators, even if the civic need for their critical perspectives increases. A valuation gap appears, whereby the demand for market-relevant critical and technical skills is often more immediate and highly compensated, which can skew educational priorities away from other objectives, almost splitting the triangle into two halves. While Ernest (1991, 2019b) already observed the marginalisation of progressive and public educators in 20th-century British mathematics education, we believe it to be an observation that can be feasibly transferred to other capitalist societies, and that it continues to be relevant in the 21st century.[8] This leads us to socio-economic tension 1:

> **The socio-economic valuation gap:** A socio-economic tension arising from a complex market dynamic. While the normative need for critical mathematical education grows to hold mathematisation accountable, and critical thinking and problem-solving are also important for industry, the dominant economic incentive structures and career pathways still appear to reward technical, optimising expertise disproportionately. This creates a valuation gap: skills essential for civic accountability and long-term sustainability are needed but often less directly or highly valued in the market than immediate technical proficiency.

The economic pressures related to the valuation effect put a complementary perspective on the curriculum tensions from the introduction. When education prioritises market-ready skills, it can leave little room for students to engage critically with sustainability issues. This is particularly problematic, as Chiodo and Bursill-Hall (2019) argue that only a multidimensional engagement with ethics in mathematics leads to success. This effect provides the mechanism for this reluctance: when educators perceive that the market overwhelmingly values technical skill, they may deprioritise ethical discussions, signalling to students that such concerns are secondary to professional success; such a concern initiated the introduction of the ESCT in the first place (Müller et al., 2025). Students who notice that their teachers and lecturers are reluctant to speak about ethics and sustainability, or feel pressured not to incorporate such concerns into mathematics exercises, may themselves become reluctant to consider such concerns.

---

[8] In this context, we wish to also point to Marktanner (2023) who studied the role that public sector consultants played in the 1990 school reforms in the German state of North-Rhine Westphalia.

## Mathematical implementations in capitalist societies

In addition to creating and structuring narratives, the increasing mathematisation involves constructing new and restructuring old markets. As Grötschel et al. (2010) argue, mathematics has solidified its position as an important new production factor, driving efficiency, setting organisational principles, and enabling new forms of economic utilisation. In particular, current trends in AI and platform economies exemplify how mathematics, capitalism, and geopolitical competition intersect. The rapid advancements in artificial intelligence have ignited narratives of a geopolitical "AI race," as nations compete for strategic advantage (Cave & Ó hÉigeartaigh, 2018; Ó hÉigeartaigh, 2025). However, this competition is not merely technological; the different national models of capitalism profoundly shape it.

Our analysis will proceed in two parts. First, we will consider a "varieties of capitalism" approach to highlight some basic tendencies. However, as this approach was developed to compare nations and does not fully capture the transnational dynamics of Big Tech, we will then complement it with a perspective on variegated capitalism. Whilst this latter framework was developed partly as a critique of the state-centric focus of the literature on varieties of capitalism, we use the two in combination: technology firms strategically exploit the institutional differences identified by varieties of capitalism, only to then supersede these, creating new hierarchical forms of economic action.

## Mathematisation meets varieties of capitalism: A story of different implementations

Hall and Soskice (2009) highlight a crucial distinction between liberal market economies (LMEs), such as the United States (US), characterised by deregulation, fluid labour markets, and reliance on stock markets for finance, and coordinated market economies (CMEs), such as Germany, which rely more on long-term strategic cooperation between firms, strong vocational training, and incremental innovation, as well as other mixed forms known as mixed market economies (MMEs). The current AI landscape — dominated by US-based firms like OpenAI (ChatGPT), Google (Gemini), and Anthropic (Claude) — illustrates that the institutional structure of LMEs possesses distinct advantages in fostering radical mathematically- and technologically-driven innovation. This advantage is particularly pronounced in the United States, where the supply-side flexibility of an LME aligns with a national economic strategy driven by aggregate domestic demand, heavily reliant on credit-fueled consumption (Reisenbichler & Wiedemann, 2022). This model is predicated on systemic financialisation and reliant on the wealth effects generated by asset appreciation in housing and equity markets (cf. Krippner, 2011; Hall & Soskice, 2009) — assets which have themselves become increasingly securitised and abstracted (Akyıldırım & Soner, 2014; Vasicek, 2016). Today's financial markets would be unthinkable without advanced mathematics, yet flawed models can have systemic consequences, as the "formula that killed Wall Street" illustrated during the 2008 crisis (MacKenzie & Spears, 2014).

Within this macroeconomic structure, capital constantly seeks new, high-yield asset classes. The technology sector, and most recently AI, has become the primary outlet for this speculative capital (Floridi, 2024). The US economy is thus institutionally predisposed to support the high-risk environment of the AI industry. It features a widespread venture capital infrastructure willing to make large, speculative bets and fluid labour markets that allow companies to scale rapidly. Furthermore, such LMEs encourage winner-take-all competition, characterised by weaker labour protections and a preference for competitive market strategies over inter-firm collaboration (Hall & Soskice, 2009). This ingrained competitive ethos is readily apparent in the strategies of leading AI companies (Crawford, 2021; Hao, 2025; Olson, 2024), the structure of Big Tech's platform economies (Srnicek, 2017), and the resulting dynamics of LME-driven forms of surveillance capitalism (Zuboff, 2019).

Liberal market economies demonstrate particular strengths in the breakthrough innovation needed to establish and maintain emerging AI markets, yet this advantage does not signal worldwide adoption of their economic framework. Mathematisation is fundamentally agnostic, and it can serve as a versatile instrument that countries can use to amplify their existing institutional strengths. Rather than homogenising global economic systems, mathematisation may actually maintain existing distinctions when nations deploy it to sharpen their specialised capabilities. For example, while AI tools are used for disruption in Silicon Valley (Geiger, 2020), they are applied differently in the German Mittelstand of high-end manufacturing. In the latter, machine learning helps to analyse complex production lines, refine manufacturing, and improve product quality (e.g., Haarmeier, 2021). However, the interdependence between advanced mathematics and the economy goes far beyond AI and finance, and it is becoming increasingly ingrained in all forms of capitalism. For example, even industries like fast fashion, a prominent case study of unsustainable economics, rely heavily on mathematics to optimise supply chains, design and produce clothes, and maximise consumption (cf. Magri & Ciarletta, 2023; Mair, 2024; Thomas, 2019). The two-fold impact of mathematics is visualised in Figure 3.

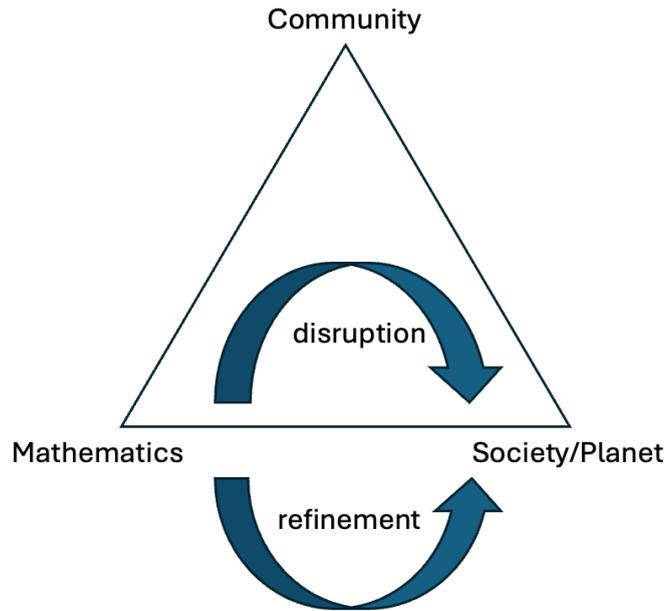

Figure 3: Mathematics as a disruptive and refining force

These differences in capitalist economies are not only connected to different perspectives on mathematics, but they can also be correlated with different educational approaches (Hall & Soskice, 2009). Regarding sustainability, several European countries have embedded strong educational goals into their national curricula (Jucker & Mathar, 2015), while the US is increasingly falling behind (DeMio & James, 2025; Klingebiel & Sumner, 2025). This may create tensions within the mathematical education community, whereby some see a need for more radical, critical positions advocating transformation, while others can be content with slower reformative processes (cf. Müller, 2024; as well as the differences in education about, for, and as sustainability presented in Müller et al., 2025). In particular, liberal market economies such as the US typically may be associated with increased pressure on community concerns that may not always be equally shared across other national education systems, as well as those from the Global South, who critique Western-centric development programmes and forms of industrialisation (cf. Müller & Chiodo, 2025). Perhaps most notably, such differences can also be observed in AI regulation, where the European risk-based approach meets more decentralised, less encompassing approaches in the US (Chun et al., 2024).

However, despite the underlying differences presented by varieties of capitalism and national economic strategies, and while economic forces may pull mathematics and its education into a more utilitarian direction, the training of mathematicians has been rather homogeneous, with little to no education in ethics and sustainability (Müller et al., 2022). As also observed earlier by P. J. Davis (1987), mathematicians still shy away from asking "embarrassing" why-questions about the process of mathematisation. As Ernest (2019a, p. 2) notes, "It would never occur to most mathematicians to question the value of mathematics, let alone to question its ethics; to ask if mathematics can possibly be harmful or damaging to individuals or society." Because of this, Thompson (2023) warns against being stuck in "model land," so as not to confuse the abstract elegance of a mathematical model with a more messy reality. To successfully train mathematics students in ethics and sustainability, one has to adjust such teaching to fit local circumstances (Müller & Chiodo, 2025). This leads us to socio-economic tension 2:

> **The multifaceted implementation of mathematics:** Just like the history of mathematics can tell both stories of homogenisation and variation in mathematical practice, different capitalist societies pull the economic implementation of mathematics into different directions, despite the subject showing socio-cultural homogenising forces. In the story of mathematics, neither homogenisation nor heterogenisation are a clear economic winner.

As our analysis has shown, any such localised education also includes the need to be aware of local socio-economic conditions and legal regimes; these will also influence how one navigates the practical tensions. For example, "Psychology 101: How to survive as a mathematician at work" (Cambridge University Ethics in Mathematics Society, 2024) means something different in a liberal market economy with little labour protection compared to a coordinated market economy with a stronger social safety net and well-maintained trade unions. This refines the "location effect" described by Müller et al. (2025) as even positions close together (e.g., two educators near the community vertex) may speak the same "language" but can still misunderstand each other because they experience different pressures, and thus understand what is said differently.

## Mathematisation meets variegated capitalism: Hierarchical dependencies between the Global North and South

While the varieties of capitalism framework of the previous section was able to hint at how mathematisation unfolds differently within core economies of Western Europe and the United States, this analysis remains incomplete, given the goal to challenge Eurocentric storytelling. Hall and Soskice (2009) treat nations as discrete and comparable units. This horizontal analysis hides the hierarchical dependencies that structure global capitalism (Peck & Theodore, 2007). To better understand mathematisation's full impact, we must then also engage with the perspective of variegated capitalism, which aims to analyse core-periphery power asymmetries and arose largely as a critique of the literature on varieties of capitalism. Where the literature in the previous section compared nations horizontally (as if comparing different systems on a level plane), research into variegated capitalism emphasises the vertical, hierarchical, and occasionally exploitative relationships between them. It recognises that global capitalism is not a level playing field; rather, powerful actors (nations and transnational corporations) actively leverage institutional differences and geographical distance for competitive advantage. This perspective aligns well with lessons from Ethnomathematics (about colonial practices in the history of mathematics) as it reveals that the clean and immaterial mathematisation celebrated in Silicon Valley depends fundamentally on hidden infrastructures of human labour and resource extraction concentrated in the Global South (as similarly argued by Crawford, 2021). Far from democratising opportunity, the mathematisation of the world is also deepening existing inequalities through new dependencies and renewed forms of exploitation.

Behind the seamless user interfaces of AI systems celebrated as pinnacles of LME innovation in the previous section lies a less visible reality. These supposedly autonomous technologies depend on invisible human workers across the Global South who label vast amounts of data for training and validation purposes (Perrigo, 2023). Such ghost work, as Gray and Suri (2019) call it, is often obscured through complex subcontracting chains and workers being paid low wages. These ghost workers may face psychological harm from prolonged exposure to disturbing content during labelling and moderation tasks, yet lack the employment protections or mental health support available to their counterparts in Silicon Valley (Perrigo, 2023). In both the Global North and South this outsourcing of tasks creates a new precariat of workers without access to the benefits of mathematisation, while being essential to it (cf. Standing, 2011; Crawford, 2021). Crawford (2021) goes a step further by outlining the materiality of these systems in regard to the entire planet — from human life to ecological concerns —, thereby turning her *Atlas of AI* into an atlas of how AI changes the (physical) world.

Crawford (2021) argues that rhetoric surrounding AI and similar forms of mathematical/digital technologies (e.g., "cloud computing") obscure their material dependencies. AI training runs, algorithmic calculations, and the collection of data require substantive infrastructure. This infrastructure, and the processes surrounding it, echo historical imbalances in resource relations, including both physical raw materials (Regilme, 2024; Crawford, 2021) and data (Luo, 2025). Just like the environmental costs remain asymmetric, AI and other algorithms regularly perform worse on data related to people from the Global South (e.g., Noble, 2018).

Viewed through the lens of variegated capitalism, the mathematisation of the world is neither a purely democratic nor a symmetric force. It does not homogenise the economic circumstances, nor does it align the living standards: those developing the mathematics and building sophisticated new technologies on it gain substantially more than those who provide the raw materials in its support. This asymmetry between how mathematical technologies such as AI materialise, can imply that certain socio-planetary threats may be easily overlooked (visualised in Figure 4 as a cloud hiding the Society/Planet vertex).

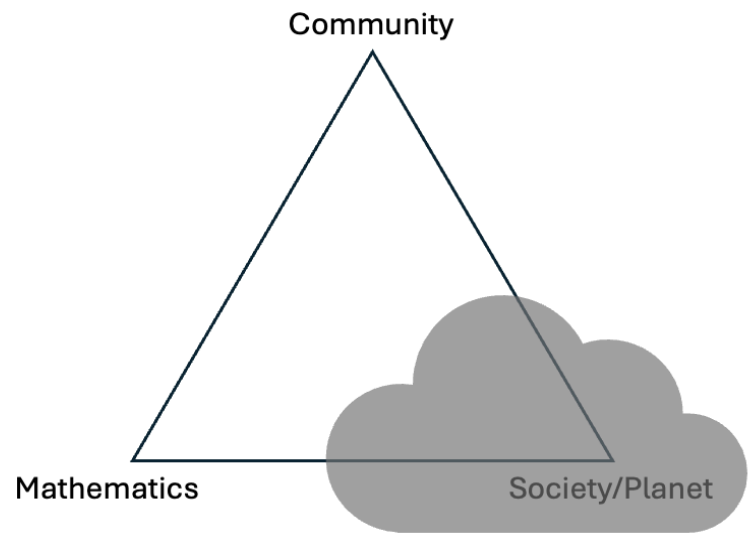

Figure 4: Hidden socio-planetary concerns

The lessons for educators are immense: educating current and future mathematics students about asymmetries between the Global North and Global South is not just a matter of ethnomathematics and past colonial histories; it is also a question of educating them about current unsustainable trends where mathematics is a substantial part of the problem. It brings to light the hidden nature of mathematics that Ernest (2019a) observed, and that societies are ill-prepared to make it visible. Decolonising a mathematical curriculum is thus not just a matter of properly contextualising the subject's history and discussing the practitioners' privilege (Borovik, 2023); it also includes epistemic decolonisation (Müller & Chiodo, 2025), such as the styles of extractive thinking enabling the creation of socially unjust data supply chains. This leads us to tension 3:

> **The material opaqueness of mathematics:** From the perspective of economics, the increased power of mathematics affects the world in many material ways, but the negative social consequences of it are often outsourced to where those building it do not see it.

Finding the correct framing of sustainability contexts is thus not only a matter of matching it with the students' expectations and their self-efficacy (as detailed by Meyer et al., 2025). The material opaqueness of mathematics is a direct pedagogical obstacle here: it becomes difficult to frame a sustainability context compellingly when its most significant negative impacts (such as resource extraction or ghost work) are invisible to both the educator and the students. This is compounded by questions about the control and power the communities in which the students live have in regard to the formatting powers of mathematics. Are they driving the innovation, utilising it, or are they even left out of their benefits completely?

Finally, we can see that mathematisation is one of the primary logics through which variegated capitalism works, and why the literature on varieties of capitalism needed to be expanded. The mathematisation of markets and economic decision-making has a lever-effect: comparatively few resources can be deployed to affect a large number of people (cf. Chiodo & Müller, 2025). McLuhan (1964) suggested that electronic media shape society more through their characteristics than through the content they disseminate. This perspective has received much criticism, but there is a deeper point about mathematisation to be made. Without mathematics, neither the electronic media described by McLuhan (1964) nor today's digital technologies would exist, nor would much of the content that gets disseminated via today's technologies. If technology is the central nervous system (CNS) of society, then mathematics enables its construction and the sending of messages. Mathematics shapes the technological CNS and what societies can do with it.

Having assembled our foundational framework through a combination of "imagined futures" (Beckert, 2016), varieties of capitalism (Hall & Soskice, 2009), and variegated capitalism (Peck & Theodore, 2007), we now wish to look at some concrete consequences of mathematising the world. What does the increased operationalisation of mathematics mean for people and institutions, including mathematics educators?

# Part III: Consequences of mathematising the world

> *"Algorithms of various kinds hold the world together. Financial transactions, dating, advertising, news circulation, work organization, policing, tasks, music discovery, hiring processes, customer relations — all are to a large extent delegated to non human agents embedded in digital infrastructure."*
>
> <div align="right">Airoldi (2022, p. 1)</div>

The transformation of society through mathematical logic extends beyond economic structures to reshape fundamental relationships between power, accountability, and human agency. As mathematical systems acquire increasing authority to define problems, frame solutions, and allocate resources, they create new forms of governance that operate through algorithmic decision-making rather than democratic deliberation. This shift generates a paradox: whilst mathematical tools promise objective, efficient solutions to complex problems, their opacity and complexity enable new mechanisms for avoiding responsibility. The technical sophistication that grants these systems their authority simultaneously complicates oversight, raising urgent questions about democratic oversight in an increasingly mathematised world.

At the individual level, this mathematical reorganisation of society reconfigures human experience itself. Algorithmic systems increasingly mediate how people understand themselves, make choices, and imagine their futures. The predictive capabilities of these systems — their ability to model behaviour, anticipate preferences, and optimise outcomes — create subtle but profound constraints on human agency. This process unfolds against a fragmenting international order, where the erosion of multilateral cooperation undermines precisely those collaborative frameworks needed to govern mathematical technologies and address their global impacts. Together, these dynamics reveal how mathematisation reshapes not just how decisions are made, but who makes them and what possibilities remain open for human flourishing.

## Between power and accountability

In the context of increasing mathematisation, Skovsmose (2021) identifies a threefold formatting power of mathematics in the public sphere. This refers to the capacity of mathematics not just to describe the world, but to actively shape and structure reality itself. This power manifests in three ways. First, mathematics can constitute a crisis (for example, complex financial instruments generating systemic risk). Second, it can picture a crisis (such as epidemiological modelling). Third, it can shape a crisis by defining the available response options: how a problem is mathematically framed influences which solutions are considered feasible, thereby exerting considerable influence on political decision-making. Regarding the ensuing climate crisis, Beckert (2024) warns that prioritising short-term profits may have severely compromised long-term sustainability, a process facilitated by quantitative models that discount future risks and well-being.

In this intertwined setting of increasingly mathematised capitalist societies, Big Tech's influence extends beyond markets into research and policy. Abdalla and Abdalla (2021) raise concerns regarding corporate influence over the AI safety research agenda, drawing parallels with historical examples of industry-funded research, including Big Tobacco. By shaping public discourse, technology companies can deflect scrutiny. For example, AI companies benefit significantly when public and policy attention focuses on preventing speculative, long-term dangers of artificial general intelligence (AGI), rather than addressing immediate sociopolitical issues arising from existing AI systems. Focusing on the existential risk dangers of AGI creates a myth of future performance surrounding these systems, which then incentivises actors to buy and use AI today (cf. Hao, 2025; Olson, 2024).

The complexity and opacity of mathematically- and algorithmically-driven systems create significant accountability challenges. While these systems are becoming increasingly better on a technological level, they can lead to the demathematisation of the population, whereby people learn to do less mathematics because the technological systems surrounding them can do more.

> "For the user of technology it becomes more important to, first of all, simply trust the black box and then, to know when and how to use it — for whatever purpose. If it turns out to be inefficient, ineffective, erroneous or disastrous, nobody can be blamed — it was the technology's fault. From this perspective, demathematisation reduces the feeling for, and the acceptance of, responsibility: a car's antilock braking system is taking as a licence for driving fast; formal assessment of personal creditworthiness ensures that a loan is not approved to the wrong people. When faults occur, then this is a request for a better technology — technology that is designed to be foolproof." (Jablonka & Gellert, 2007; building on Keitel, 1989; Keitel et al., 1993).

Here, too, mathematisation mirrors general organisational principles of capitalist societies. Porter (2020) argues that the increasing quantification of the socio-economic sphere is a deliberate political choice, often coming out of a weakness of those in power and in order to establish trust through a seemingly neutral language of numbers. Along with such "mechanics of objectivity" and "technologies of distance" comes rigidisation, whereby an organisation becomes more rigid as it has to now adhere to certain (fixed) standards, rules, and processes to enable the use of mathematics (cf. Porter, 2020). The adherence to these, however, can quickly turn into a problem of accountability, allowing institutions and individuals to evade responsibility.

To capture this, D. Davis (2024) introduces the concept of "accountability sinks" to describe how industries and political systems evade accountability by hiding behind complex processes and pre-defined rules, which often fail to allow for exceptions and human agency when needed. An "accountability sink" effectively absorbs responsibility, distributing it across a complex system until no single individual can be held liable for failures. As an example from education, consider an algorithmic university admission system that, unfortunately, makes biased decisions. The admissions office may blame the software; the software developer blames the flawed training data; and the institution claims adherence to standardised, objective protocols. If there are no rules in place that empower the applicant to legally challenge the decision, the complexity of the mathematical system becomes a shield, leaving them without recourse.

But while accountability sinks are already an issue today, Chiodo et al. (2025) go a step further. They warn that we may be approaching a time when governments "hand over the keys to the city" to private, non-elected groups of experts and industries, including those mathematically trained. This transfer of authority may bypass democratic deliberation and oversight, raising concerns about the future of democratic control in increasingly mathematised societies. The concentration of socio-economic power in a select group of mathematically trained experts creates new political tensions. In particular, it stands in sharp contrast to the potential "collateral damage" that Ernest (2019a) identifies regarding the learning of mathematics. Ernest's ethical audit of the subject highlights that for many individuals, mathematics is not a tool for empowerment, but rather "cold", "hard", and "unforgiving", and a deeply alienating experience. We want to suggest here that Ernest's analysis is not merely an individual psychological concern (i.e., maths anxiety) or the beginning of a burgeoning socio-economic divide, but that it can threaten the very functioning of our democratic societies.

As society becomes intensely mathematised, its overvaluation and filter function do more than create individual failures. It risks manufacturing a new, politically defined stratum of citizens who feel (and often are) left behind by the very subject that they already despised at school. This dynamic would threaten to make mathematics itself, including but not limited to AI and algorithms, a central signifier in and tool for populist narratives. Such a process mirrors what Di Nucci (2020) identifies as the "control paradox": the act of delegating control to complex (mathematical) systems, and the experts controlling these, makes citizens feel disempowered and neglected, which in turn motivates a populist backlash. It pits a remote, often unaccountable expert class against a public that not only mistrusts but has been actively alienated by the very language in which this new form of power is expressed (see also Berman & Hirschman, 2018). Skovsmose (1998) similarly notes that the increasing mathematisation of the world creates a "risk society", potentially happening outside of democratic control, and threatening to leave those behind who lack the socio-economic power and mathematical knowledge to question the (often hidden) formatting power of mathematics.

As Chiodo et al.'s (2025) "handover of keys" happens, it may have particularly damaging impacts in high-stakes domains such as criminal justice, which then meets an already mistrusting group of citizens. For example, mathematically-driven predictive policing systems often display systemic biases, leading to overpolicing of already disadvantaged areas (J. Davis et al., 2022; Egbert & Mann, 2021). Similarly, recidivism algorithms used to predict whether defendants will reoffend have exhibited racial bias. These systems have sparked intense legal debates over procedural due process (Washington, 2018). Critiques that focus only on the quantification involved miss a deeper concern: the mathematisation of justice presumes that justice cannot just be optimised using mathematics, but that justice itself is a natural phenomenon which can be properly captured by the language used to control and exert power over nature, i.e., mathematics. It is this underlying logic of optimisation and control that constitutes the primary ethical problem, existing even before a single data point is collected. What would a person who already despises mathematics make of this?

This increasing reliance on mathematics raises a further fundamental question: when human behaviour appears as a machine-learning outlier, can we always remedy the discrepancy through more data and more complex models, or are humans inherently too complex and varied ever to become fully capturable by mathematics (Marichal, 2025, p. 7)? As Marichal argues, the technology sector appears to assume the former, driving a relentless pursuit of data collection and modelling. In contrast, sociologists and AI safety scholars often appear to assume the latter. Perez (2020), for example, highlights how data on women is often missing in many domains, even in everyday situations like public transport planning, which leads to systems designed by default for men. While significant gaps exist in the data we collect, the question becomes: Does it matter what the answer is? As LLMs are hitting diminishing returns in scaling laws (Chen et al., 2025), it may very well be the case that both perspectives lead to imperfect models, unable to capture the behaviour of outliers properly.

The fact that these imperfect[9] models and forms of mathematisation are still good enough to be used increasingly becomes a problem in civic life. For example, algorithmically curated information streams and bots shape civic discourse and political awareness (Chang et al., 2021; Gensing, 2019). As algorithms and other forms of mathematics increasingly mediate social interactions and economic decisions, individuals find themselves categorised (Marichal, 2025; O'Neil, 2016). And as democratic societies have to deploy increasingly strong rules and regulations to keep up with the changes (Deneen, 2019), and some are even actively considering using AI to make decisions on violence and war (Erskine & Miller, 2024), the risk of unaccountable automated forms of governance may further increase in the future.

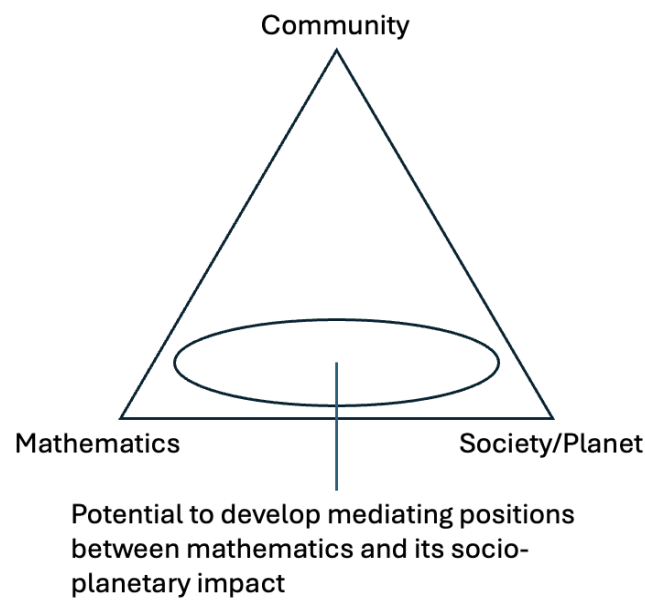

Figure 5: The need for mediation between mathematical concerns and socio-planetary impact

This presents the mathematics education community with a challenge, whereby critical perspectives along the mathematics-society/planet edge have been less developed than those along the community-society/planet edge (cf. Müller et al., 2025). However, the "location effect" described by Müller et al. (2025) suggests that it is precisely the former discourses which may be best suited to mediate between desires to mathematicise the world and the needs of today's democracies (Figure 5), bringing us to socio-economic tension 4:

> **The economic power versus social unaccountability gap:** The socio-economic tension that exists between the escalating economic and formatting power granted by complex mathematics, and the persistent democratic necessity of holding that power to account. This conflict manifests as unaccountability gaps, where opacity and complexity are leveraged to create systemic structures (such as accountability sinks) that shield economic actors from social and democratic oversight.

---

[9] We do not want to suggest here that every imperfect model is a bad model. There are many case studies where introducing algorithms have actually improved the situation. Berman & Hirschman (2018, pp. 264 - 265) discuss such an example of algorithmic bail reform. We do, however, see a problem when imperfections become systemic.

This unaccountability gap potentially fuels polarisation. On the one hand, routinisation of work is known to be able to reduce ethical awareness and actions (Kelman, 1973). On the other hand, we see a mechanism of alienation: when students (particularly those from marginalised groups) experience the unaccountable power of a biased algorithm as a neutral, mathematical fact, it can foster a sense of learned helplessness. As such, when they confront the unfairness of these seemingly neutral yet unaccountable systems, it may affect them negatively. Following D. Davis's (2024) exposition of what happens when people encounter such systems, we suggest that two extremes are possible: in the worst cases, students enter into a rejection of mathematics' utility (as they do not experience its benefits) or into an acceptance of its authority (as they cannot properly challenge it). As Ernest (2019, p. 2) notes, "for a significant minority of the public at the very least, mathematics is seen to be cold, hard, unforgiving, masculine, meaningless, joyless, rejecting or frightening," and their encounter with mathematically induced accountability sinks may only make such matters worse.

## The challenge to human identity

Processes of datafication and mathematisation create new forms of identity and shape our understanding of freedom. Thus, Marichal (2025, p. 14) argues that "algorithms frustrate the principle of liberal autonomy in that they interfere with our ends by predicting and anticipating our desires, possibly before we even know them ourselves", and that this preemptive shaping of behaviour circumvents conscious decision-making, suggesting a future where choices are optimised for us rather than made by us. This brings us to the central paradox of mathematisation. Beckert (2016) suggests that mathematical models enable economic actors to expand and structure possibilities and make decisions. However, when applied outside of economic transactions and to the private lives of citizens, mathematisation can have the opposite effect. The optimisation-driven processes described by Marichal (2025) and the general mathematisation processes described in this paper enable the enclosure of imagination.[10] The potential long-term consequences, thus, may very well turn Beckert's (2016) observations upside down by reducing the human capabilities to imagine and think in alternative futures.

Rather than assisting in and legitimising human decision-making, mathematically- and information-driven technologies may increasingly restrict human control (Chin, 2025). Zuboff (2019) similarly argues that the "behavioural surplus" (i.e., the availability of masses of user data) increasingly turns users into something to be predicted and monetised for marketing purposes. As Marichal (2025) further notes, the better humans conform to these predictions, the better it is for the firms relying on them. Thus, in particular, being a human "outlier" is contrary to the platform business and other statistically-driven market interactions — in other words, the optimisation principles underlying these economic transactions begin to put pressure on both individual agency and the ability to imagine (global) sustainable alternatives (Figure 6). Consequently, Marichal argues that optimisation does not just describe reality through quantification and change this reality, but that optimisation has a disciplining function on the individual. In other words, an increasingly strong focus on mathematisation threatens to narrow the range of imaginable futures to only those which can be made legible through mathematics itself. The logic that frees capital and enables economic decision-making simultaneously restricts human agency in other circumstances.

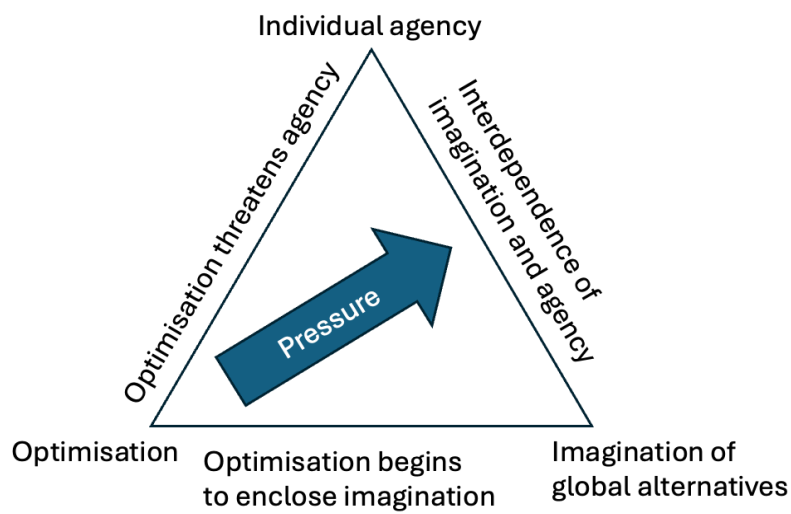

Figure 6: The impact of optimisation

---

[10] The term is inspired by Buschert and Kujundzic (1997).

Higher mathematics increasingly affects potentially life-defining decisions (O'Neil, 2016). When these decisions become automated and opaque, individuals often lack the means to understand or contest the outcomes (Chiodo et al., 2024). These "weapons of math destruction" (O'Neil, 2016) may then systematically disadvantage certain groups and exacerbate existing inequalities. To capture these subtle forms of power, Cheney-Lippold (2011) reconceptualises Foucault's perspectives on biopolitics into "soft biopolitics" to describe the indirect power and identity changes occurring through algorithms.

While mathematical models pervade the public sphere, they are often accompanied by claims of neutrality and objectivity. However, technologies are not neutral, as Winner (1980) famously argued: Technological artefacts have politics. Müller and Chiodo (2023) argue that this principle applies equally to mathematical artefacts. Among mathematics, optimisation has become the dominant mathematical paradigm shaping the public sphere. Algorithms typically maximise a specific objective function, whether profit, efficiency, or engagement. However, optimisation inherently involves trade-offs and value judgments, as deciding what to optimise and what constraints to impose is a political act (Chiodo & Müller, 2025; Marichal, 2025; Müller & Chiodo, 2023)

What we see here is that while different varieties of capitalism and national economic strategies do exist (Baccaro et al., 2022; Hall & Soskice, 2009; Pontusson et al., 2021), and different regulatory regimes have been developed to deal with the threats from AI (Chiodo et al., 2024), there are still strong commonalities and homogenising forces — particularly, regarding the effects on the identity of individuals and on the need for civic education. Mathematising humans can challenge their identity at a fundamental level, and as Carvalho et al. (2005) argue, it is not enough to consider their rights and duties; citizenship and being human are more than that. Critically, however, the effects of mathematising human identities demand teachers with a strong ethical awareness, as those with higher levels are better equipped to enact meaningful change in educational institutions and advocate for their students (Rycroft-Smith et al., 2024). This brings us to socio-economic tension 5:

> **The enclosure of imagination:** The increasing mathematisation of the socio-economic world provides an increasing number of possibilities for deploying mathematics; however, this increased deployment does not necessarily create a more imaginative future. Instead, it threatens to restrict our thinking in sustainable alternatives.

The enclosure of imagination directly threatens the capacities to deal with complex ambiguities. This is because an optimisation paradigm relentlessly trains individuals to seek a single, quantifiable, correct solution, thereby atrophying the very skills — such as managing trade-offs, embracing ambiguity, and imagining non-linear solutions — that are essential for addressing complex sustainability problems (cf. Chiodo & Müller, 2025; Müller & Chiodo, 2023). In particular, the optimisation paradigm deployed in capitalist societies contrasts sharply with the imaginative practices necessary in teaching ethical and sustainable mathematics. The enclosure described here may partially explain why the field struggles with what Makramalla et al. (2025) call the need to imagine new forms of sustainable mathematics education. In particular, the present forms of mathematisation threaten to confine the aims of education to the quadrant of "economic productivity", rather than "stewardship of the planet", "reciprocity & collectivism", and "activism" (Makramalla et al., 2025, p. 547) — quadrants with a stronger focus on sustainable mathematics education.

# The decline of multilateralism and its challenge to critical paradigms

## A critical pragmatic appraisal of the situation

While the mathematisation of (all) parts of society continues, the unique system of multilateralism created post-World War 2 is deteriorating, and a realist perspective would be to get accustomed to alternative, potentially more fragile forms of international collaboration (Bekkevold, 2025). However, as multilateral efforts are particularly important to regulating AI and other forms of mathematics (Natorski, 2025; Pfeiffer, 2022), and institutions like the United Nations play a significant role in setting a common language on sustainability (Biermann et al., 2022; Biermann et al., 2017), their decline can have direct effects on ethical and sustainable mathematical practice and its teaching.

One can observe that the perspectives from critical mathematics education, which underlie most current approaches to sustainable and socially-just mathematics education, have pedagogical perspectives and problem-solving approaches that are strongly aligned with multilateral and pluralistic decision-making (Table 1). Thus, the decline of multilateralism presents a particular challenge to public and progressive educators, as many of their goals, pedagogical approaches, normative grounding, and discursive practices may be called into question. This shift in the discursive landscape and in the grand strategies of international politics meets an already eclectic mix of tensions between critical and absolutist stances on mathematics. We suggest here that this may negatively affect normative positions near the community-society/planet edge in the ESCT. As the geopolitical landscape fractures and nations retreat into more realist, interest-driven postures (Bekkevold, 2025), the societal support for these critiques potentially erodes. Educators focused on the Community-Society/Planet nexus (Figure 7) may find their perspectives marginalised in increasingly nationalistic or purely utilitarian environments, where arguments grounded in global solidarity or universal rights lose traction against immediate economic or security interests.

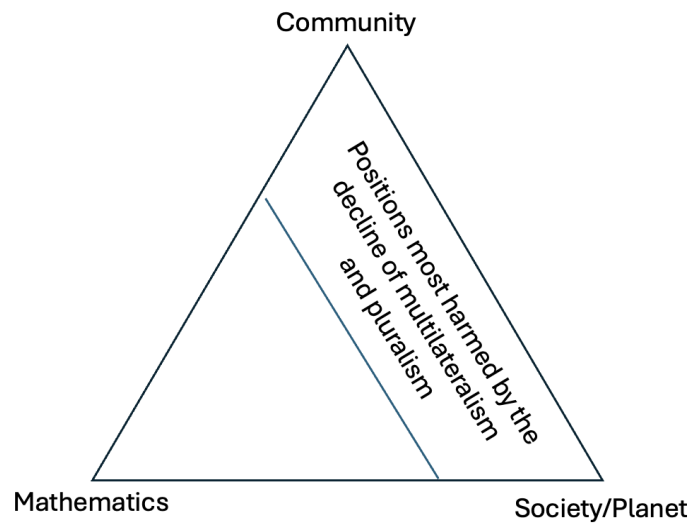

Figure 7: Potential effects of declining multilateralism

Associated with the decline in multilateral approaches, we see an imperative for critical pragmatism and for deeper localised perspectives on ethical and sustainable mathematics (cf. Müller, 2025a; Müller & Chiodo, 2025). But while Müller and Chiodo (2025) argue that this localisation need of ethics and sustainability comes from having to connect with the existential and cultural concerns of students, teachers, and professional mathematicians, there is a second dimension to this: with increased international fragmentation, global frameworks are increasingly difficult to envision, implement, and enforce. As reliance on overarching frameworks and laws becomes increasingly difficult, the importance of individual ethical and sustainable awareness grows. While this highlights the importance of sociopolitical scholarship on mathematics education (e.g., the sociopolitical turn described by Gutiérrez, 2013), it also brings to light some of its limitations, as it does not put a strong emphasis on the necessity of individual actions as later, more ethics-focused scholarship (Müller, 2025b). This leads us into socio-economic tension 6:

> **The erosion of multilateralism undermines the critical challenger:** There exists a socio-political tension between the global, collaborative ethos inherent in critical mathematics education (which is necessary to address transnational problems) and the

> opposing geopolitical reality of eroding multilateralism. This conflict may create a dynamic that strengthens utilitarian-absolutist viewpoints while threatening to marginalise critical perspectives along the Community-Society/Planet edge, potentially forcing practitioners into a difficult position. The well-articulated need for localised ethical action can then look like a retreat, as many of the challenges critical educators face continue to have a global dimension.

The erosion of multilateralism transforms this pedagogical question from a theoretical one into a pragmatic one. When global frameworks (like the SDGs) lose their unifying force, an educator's reliance on them can feel disconnected from the more nationalist or realist rhetoric students are exposed to on social media and elsewhere. This geopolitical shift adds a practical reason to be cautious, compounding existing critiques of such frameworks. For example, while Müller and Chiodo (2025) argue that the Sustainable Development Goals (SDGs) are a good way to design exercises, they also note a critical downside: they follow the logic of Western development programmes. The decline in multilateralism, therefore, also suggests it is beneficial to focus on other forms of embedding ethics and sustainability for purely practical reasons. This creates a broader tension for integration across different national curricula: how much practical inspiration can, and should, be derived from other national contexts when designing mathematics tasks focused on ethics and sustainability? Where do global and local needs meet and contradict each other?

## A philosophical appraisal of the situation

The socio-economic tensions described throughout this paper, particularly the challenge posed by the decline of multilateralism, are deeply connected to a fundamental philosophical divide regarding the nature of reason and governance. This divide separates the technocratic tendencies inherent in much of the mathematisation process from the deliberative ideals championed by critical mathematics education. The philosophical underpinning of the technocratic approach can be traced back to what the early Frankfurt School thinkers Adorno and Horkheimer (1997, [1944]) called "instrumental reason" in their work on the dialectic of enlightenment. Instrumental reason describes a narrow, calculating logic that focuses solely on finding the most efficient means to reach predetermined goals. While this form of reason excels at prediction and control, it systematically sidesteps critical reflection on whether those goals are morally or socially worthwhile.

We have seen that mathematisation, when driven by this instrumental reason, tends to create governance structures that are both technocratic and hierarchical. When mathematics is viewed through an absolutist lens — as purely objective, universal, and neutral — we risk conflating mathematical optimisation with rationality as a whole. This dynamic echoes Porter's (2020) argument that claims of "mechanical objectivity" are often deployed politically to replace contentious debate with seemingly impersonal administration. This transforms complex socio-economic challenges from political and ethical questions into technical problems that seem best solved by experts; thus, leading to the "handover of keys to the city" that Chiodo et al. (2025) worry about. Within this stance, prioritising mathematical expertise over democratic participation appears not only as effective and efficient, but may also be seen as the only rational path forward.

Adorno and Horkheimer argue that the Enlightenment project, when dominated by instrumental reason, results in a paradox: the very reason meant to free people leads to new forms of domination; what was meant to be an emancipatory force becomes the very mechanism that exerts control over people and nature. The mathematisation of the socio-economic sphere that we described in the previous sections can be understood as the 21st-century manifestation of this dialectic. When mathematical systems are governed technocratically, i.e., with a focus on utilitarian and other instrumental goals, their complexity becomes a political mechanism of control. The opaque decision-making structures and "accountability sinks" described by D. Davis (2024) are consequences of this deeper style of thinking — one that traces back through the Enlightenment to Ancient Greek methods of geometric proof (cf. Bursill-Hall, 2002). Such a mode of governance is particularly compatible with neoliberal ideologies, which rely heavily on the authority of supposedly neutral metrics and expert management to organise the socio-economic sphere. This again underscores why we speak of mathematisation rather than quantification: it is about a style of thinking, not just the production and use of numbers. In short, these systems exemplify the Frankfurt School's central concern: the triumph of impersonal administration over individual autonomy and democratic accountability.

But what about critical mathematics education?

Skovsmose (1998, p. 196) notes that,

> "The claim of critical education is: Fundamental undemocratic developments must be challenged in education! In fact, the development of 'critical education' can be seen as an educational attempt to provide a new foundation for education for citizenship."

As the socio-critical stance underpinning critical mathematics education traces its roots back to the Frankfurt School (Müller, 2024), it has embedded in it an explicit critique of instrumental reason. The predecessors of critical mathematics education scholars defined critical theory in opposition to logical-positivist theories and myths of instrumental neutrality, and authority in particular (Schwandt, 2010). Critical perspectives on education insist that all knowledge and human practice, including mathematics, is shaped by social interests and power relations (Gottesman, 2016). Consider, for example, the work of Skovsmose (2021), who views mathematics as a human activity possessing significant formatting power, including the ability to constitute a crisis. As a consequence, reason and all (a priori and a posteriori) formal expressions (like mathematics) have politics embedded within them.

As another example, consider Atweh and Brady's (2009) approach to "socially response-able" mathematics education, which builds on Levinas and Habermas to develop an ethical grounding of mathematics education in the need to properly communicate and respond to the needs of the Other. To counter the dominance of instrumental reason, second-generation Frankfurt School scholars like Habermas developed alternative models. Habermas (1995 [1981]) proposed a governance model based on communicative (discursive) rationality. In contrast to instrumental reason, he argued that consensus and mutual understanding, achieved through reasoned argumentation and discourse, are the morally best way to govern. While critical education draws heavily on this deliberative ideal, it also incorporates perspectives (such as those inspired by Foucault) that emphasise how power dynamics can distort discourse, necessitating constant attention to inclusion and the critique of established norms (cf. Gottesman, 2016). However, recall that even Habermas saw mathematics as largely harmless (Skovsmose, 1998, p. 195). Critical scholarship, therefore, mandates that the application of mathematics in society must be subject to deliberative, democratic scrutiny, but in its mission, the field has developed an even more critical outlook on mathematics than its philosophical forbears.[11] From the perspective of critical mathematics education, technical correctness, as judged by experts, is insufficient for legitimacy; coercion, manipulation, and exclusion are instrumental forms of reasoning that run against this vision; and mathematics itself is not a neutral form of rationality.

This brings us back to the geopolitical situation. The two philosophical stances imply different models of governance (as summarised in Table 1). The technocratic approach, reliant on instrumental reason, is strengthened by the decline of multilateralism because, in the absence of robust global democratic deliberation, the appeal of seemingly neutral, expert-driven solutions grows. Conversely, critical education scholars trace their philosophical groundwork to legitimacy derived through deliberative, democratic means. The public sphere is a cornerstone that needs to be protected and strengthened, and multilateralism represents the primary mechanism for realising this deliberative ideal in international relations. Coming back to the mathematisation of the socio-economic sphere, this means that such a process can never be legitimised through technical validation alone, but only through communication and consensus. The aim of critical scholarship is to ensure robust democratic oversight, where mathematisation can be openly debated and contested. However, it is precisely this deliberative foundation that is threatened by the decline of multilateralism.

---

[11] We wish to note that Gottesman (2016) similarly argues that critical education has undergone a long transformation which calls from "Marxist critique to poststructuralist Feminism to Critical Theories of Race."

This philosophical divide helps to situate the call for a more critical pragmatism within mathematics and its education. It is crucial to distinguish this from the instrumental reason that underpins technocratic governance. Our proposal of critical pragmatism is not aimed at serving the existing socio-economic system more efficiently. Rather, we draw on pragmatism to find a path for praxis that allows for concrete, situated ethical action in response to the critical analysis of the structural tensions previously outlined.

Critical pragmatism is not merely a call for "what works"; rather, it is a direct response to the "location effect" observed by Müller et al. (2025). This effect highlights that different positions within the Ethical and Sustainable Concerns Triangle embody fundamentally different — and often conflicting — epistemological, ontological, and ethical stances. These deep-seated philosophical divides co-create the very tensions that can pull discourses with strong mathematical concerns (which may lean toward instrumentalism) and those with strong community concerns (which may lean toward pure critique) apart. As argued in Müller and Chiodo (2025), for mathematics to become a fully democratic force, a transformation of both the practice and its education is required. Such a task is difficult when the "location effect" makes opposing philosophical positions seem mutually incomprehensible.

We wish to extend this argument here by suggesting that the educators located along the Mathematics-Community edge of the ESCT exemplify aspects of critical pragmatism. Their position requires a delicate philosophical balancing act. On the one hand, they are critical enough to reject purely instrumentalist, absolutist philosophies of mathematics. This allows them to connect meaningfully with community concerns for student agency, well-being, and the social context of learning (goals often dismissed by a purely technocratic stance). On the other hand, they remain pragmatic enough to value the robust power and utility of mathematical knowledge itself, avoiding a position of pure critique that might devalue the subject entirely. They are also critical enough to see both the moral good and bad mathematics can do, and pragmatic enough to understand that achieving good requires equipping students with strong mathematical skills, which can then be deployed to do good.

The preceding analysis detailed a deep philosophical divide between instrumental and deliberative approaches to reason and governance; this divide is synthesised in Table 1. It is important to emphasise that this table contrasts idealised philosophical traditions; it is not intended as a categorisation of individuals. As Müller et al. (2025) show, mathematicians and educators typically hold nuanced positions, resembling a "convex combination" of three core ethical concerns (mathematics, community, socio-planetary) and blending aspects of different traditions. Therefore, the table should be read not as a critique of any individual practitioner, but as a map of the ideological landscape, highlighting the poles of the philosophical conflict that structure the field., and why the decline in multilateralism may strengthen the absolutist stance, while threatening to weaken the critical stance.

|  | **Absolutist Stance** | **Critical Stance** |
|---|---|---|
| **Nature of Mathematical Knowledge** | Mathematics is objective, universal, and certain. Truths are discovered and exist independently of human consciousness or culture (Ernest, 1991a, 1991b; Müller & Chiodo, 2025; Müller, 2025b). | Mathematics is a human activity, culturally situated, and fallible (open to revision). Knowledge is constructed and embedded in social practices (D'Ambrosio, 2016; Ernest, 1991a, 1991b; Gutiérrez, 2013; Müller, 2024; Skovsmose, 2023). |
| **View on Neutrality and Values** | Mathematics is viewed as inherently neutral and value-free (Ernest, 1991a, 1991b; Müller & Chiodo, 2025; Müller, 2025b). Focus is primarily on technical accuracy, optimisation, and efficiency (Chiodo & Müller, 2025; Marichal, 2025; Müller et al., 2022). | Mathematics is inherently value-laden and political. Mathematical artefacts embody the values, interests, and biases of their creators (Chiodo & Müller, 2025; Gutiérrez, 2013; Marichal, 2025; Müller & Chiodo, 2023; Skovsmose, 2023). |
| **Pedagogical Approach** | Emphasis on certainty, procedural fluency, rule-following, and the transmission of knowledge (Ernest, 1991a, 1991b; Threlfall, 1996). History is often presented Eurocentric (cf. Hodgkin, 2005; Joseph, 2011) | Emphasis on student-centred teaching, critical thinking, comfort with ambiguity, and inquiry into diverse mathematical practices, community, and socio-planetary concerns (Carvalho et al., 2005; Gutiérrez, 2013; Skovsmose, 2021) |
| **Problem-Solving Focus** | Solving mathematically clear, often decontextualised problems where mathematics is the primary tool for the objective solution (Chiodo & Müller, 2025). Tends to prioritise optimisation of quantifiable elements (e.g., profit, efficiency, utilitarian ethics) (Chiodo & Müller, 2025; Marichal, 2025; Müller & Chiodo, 2023, 2025). | Investigating complex, ill-defined real-world problems. Analysis of how mathematics may constitute, picture, or shape the issue itself (the "formatting power") (Gutiérrez, 2013; Müller & Chiodo, 2025; Müller et al., 2025; O'Neil, 2016; Skovsmose, 2021, 2023). |

| Relationship to Power | Tends to reinforce existing power structures and control (cf. Bishop, 1990) and used for legitimising decisions while potentially hiding motivations (Carvalho et al., 2005; Porter, 2020). | Analyses how mathematics formats society and distributes power; seeks to challenge inequalities and promote social justice (Gutiérrez, 2013; Müller et al., 2025; O'Neil, 2016; Skovsmose, 2021, 2023). |
|---|---|---|
| Implied Model of Governance | May align with technocratic, hierarchical, and neoliberal governance by equating quantification with rationality. Through purely instrumental reason, it can be prone to creating opaque systems, accountability sinks, and favour expert expertise over democratic deliberation. | Advocates for pluralistic, deliberative, and democratic oversight of mathematical systems. Specifically developed to counteract absolutist perspectives, and instrumental forms of reasoning — aims at a deeper civic engagement. |

Table 1: Comparison of (idealised) absolutist and critical stances

# Conclusion

## What did we achieve?

This paper has examined how the mathematisation of the world constitutes a fundamental reorganisation of societal structures, whereby mathematical logic becomes the predominant framework for ordering and legitimising social and economic practices. Through coercive, mimetic, and normative institutional isomorphisms, this transformation extends beyond mere technical applications to reshape the cognitive devices through which futures are imagined and decisions justified. Following in the footsteps of Skovsmose's (2021) work, our examination connected the formatting power of mathematics — constituting, picturing, and shaping crises — with lessons from economic sociology to present a more refined picture. To structure this analysis, we identified and explored six key socio-economic tensions arising from this process. By focusing on mathematisation, we moved the analysis from a quantification and data-centric focus to a political-economic one, in particular questioning why mathematics is increasingly given the power to impose specific governing logics and who profits from it. In other words, we brought Chiodo et al.'s (2025) warning of "handing over the keys" to technical experts to the case study of mathematics.

First, we described a common utilitarian trajectory favouring market-oriented educational archetypes whilst marginalising critical civic perspectives. We visualised this through the Ethical and Sustainable Concerns Triangle and its fracturing along a socio-economic centre line. Subsequently, we complemented this story of homogenisation with one highlighting the heterogeneous nature of capitalist societies. By employing a varieties of capitalism framework, we highlighted that mathematisation is not always a global homogenising force; it can amplify existing institutional and national strengths within different capitalist societies. This analysis revealed a tension regarding the multifaceted implementation of mathematics, as liberal market economies such as the US currently attempt to harness radical innovation through speculative AI development, whilst coordinated market economies deploy such tools to refine their existing economies. However, this analysis proved incomplete in capturing the transnational dynamics of mathematisation. Thus, we complemented this horizontal analysis with the perspective of variegated capitalism to reveal the hierarchical, core-periphery power asymmetries at play. We argued that the supposedly immaterial mathematisation of the Global North depends fundamentally on extractive infrastructures of human labour and resources concentrated in the Global South, leading to a tension on the material opaqueness of mathematics.

Next, we discussed how mathematisation can go hand in hand with unaccountability and opacity, a dynamic that may threaten democratic governance by circumventing deliberative processes through the deployment of mathematical experts. This created a tension between economic power and social unaccountability. In this context, we discussed how human identities are affected by the increases in mathematisation through what we termed the enclosure of imagination, as the human capacity for imagining and pursuing alternative futures diminishes. We argued that these developments demand particular attention from mathematics educators, who must navigate between economic pressures for technical training and societal needs for critical literacy.

Finally, we concluded that the deterioration of multilateral frameworks compounds these challenges, undermining the global coordination essential for regulating mathematical technologies and addressing sustainability imperatives. We showed that this geopolitical fragmentation particularly affects critical paradigms in mathematics education, whose pluralistic foundations align poorly with increasingly realist international relations. This constitutes the final tension, as the erosion of multilateralism puts pressure on the challenger, i.e., critical mathematics education.

In response to these tensions, we suggested focusing on critical pragmatic approaches that acknowledge the localised nature of ethical practice whilst maintaining a reasonable commitment to broader global principles of justice and sustainability.

# Where to go from here?

The paper was positioned as an exploratory investigation, examining the potential value of connecting economic sociology with debates on ethics in mathematics education. While the framework and the six tensions identified provide an initial mapping of this complex terrain, they simultaneously reveal substantial areas requiring further investigation, spanning empirical, theoretical, and practical domains.

First, the theoretical foundations established here require both deepening and validation. The socio-economic tensions presented in this analysis are forward-looking, attempting to anticipate future challenges rather than merely describing past events. As such, a critical next step is their robust empirical validation across diverse economic and educational contexts. Concurrently, the engagement with economic sociology warrants considerable further research. We have only scratched the surface of what concepts like institutional isomorphism, imagined futures, and variegated capitalism can bring to the field. Subsequent work should deepen these perspectives and integrate additional sociological theories to enrich studies on the ethics of mathematics education.

Second, the scope of inquiry into the mathematisation of the world must be broadened. Due to spatial constraints, our analysis of variegated capitalism and the material opaqueness of mathematics focused significantly on AI. Further empirical and theoretical research is required into how mathematics "materialises" in different economic systems beyond the digital realm, and what its ensuing effects are on education. This necessitates a deeper engagement with postcolonial and decolonial frameworks to understand how mathematical practices and education either perpetuate or challenge hierarchical economic structures, and to identify the ethical questions that must be prioritised.

Third, the macro-level tensions identified demand translation into concrete pedagogical practices and methodological tools. However, the higher-order tensions presented here, complementing the concrete pedagogical tensions identified elsewhere, require more than individual responses. Like many issues of ethics in mathematics, they are also issues of policy. In particular, the enclosure of imagination presents an urgent challenge, necessitating educational approaches that actively cultivate new, divergent forms of mathematical thinking capable of envisioning sustainable alternatives.

Furthermore, research should develop and evaluate strategies that enable practitioners not only to recognise model limitations but also to navigate the "accountability sinks" related to mathematics. We outlined how these accountability sinks may extend beyond algorithms and AI, relating also to absolutist or traditionalist mathematical mindsets. We therefore suggest a closer inspection of how mathematical values can manifest in institutional structures and processes, even when explicit mathematical work is not being done.

Finally, the suggested turn towards critical pragmatism (rather than purely critical stances) amid eroding multilateralism necessitates methodological innovation. If ethical action becomes increasingly localised, researchers must continue to develop approaches that effectively connect macro-level insights with classroom realities and the lived experiences of students.